
\input psfig
\input amssym.def
\input amssym
\magnification=1100
\baselineskip = 0.25truein
\lineskiplimit = 0.01truein
\lineskip = 0.01truein
\vsize = 8.5truein
\voffset = 0.2truein
\parskip = 0.10truein
\parindent = 0.3truein
\settabs 12 \columns
\hsize = 5.4truein
\hoffset = 0.4truein

\setbox\strutbox=\hbox{%
\vrule height .708\baselineskip
depth .292\baselineskip
width 0pt}
\font\caps=cmcsc10

\def\sqr#1#2{{\vcenter{\vbox{\hrule height.#2pt
\hbox{\vrule width.#2pt height#1pt \kern#1pt
\vrule width.#2pt}
\hrule height.#2pt}}}}
\def\square{\mathchoice\sqr46\sqr46\sqr{3.1}6\sqr{2.3}4}
\def\leaderfill{\leaders\hbox to 1em{\hss.\hss}\hfill}
\font\bigtenrm=cmr10 scaled 1400
\tenrm

\centerline{\bigtenrm THE CANONICAL DECOMPOSITION}
\centerline{\bigtenrm OF ONCE-PUNCTURED TORUS BUNDLES}
\tenrm
\vskip 14pt
\centerline{MARC LACKENBY}
\vskip 18pt

\centerline{\caps 1. Introduction}
\vskip 6pt

Epstein and Penner proved that any complete, non-compact,
finite volume hyperbolic manifold admits a canonical 
decomposition into hyperbolic ideal polyhedra [2]. This is an 
important construction, particularly in dimension three, and
yet it remains mysterious. Although it can be computed in practice
for any particular hyperbolic 3-manifold, as demonstrated by the computer
program SnapPea [11], it is difficult to determine the
decomposition of an infinite family of examples. In this paper, we do
exactly that for once-punctured torus bundles over the circle
that admit a hyperbolic structure.
They have a natural ideal triangulation, first constructed by
Floyd and Hatcher [3], which we term the {\sl monodromy ideal
triangulation}. It is described in detail in \S2 of the paper.
Our main result is the following.

\noindent {\bf Theorem 1.} {\sl The canonical polyhedral decomposition
of a hyperbolic once-punctured torus bundle over the circle 
is its monodromy ideal triangulation.}

Hyperbolic structures on once-punctured torus bundles were first
studied by J{\o}rgensen. He produced an unpublished manuscript on 
quasi-Fuchsian once-punctured torus groups [5], but this did not deal with 
bundles. More recently, Akiyoshi, Sakuma, Wada and Yamashita have
developed his methods further [1], and Akiyoshi has a programme
for dealing with the bundle case.
All these approaches have been highly geometric.
Our techniques, however, are much more topological. They draw on
Gabai's concept of thin position [4] and almost normal surface
theory, due to Rubinstein [8], Thompson [10] and Stocking [9].

A once-punctured torus admits an involution, as shown in Figure 1,
which commutes with any linear homeomorphism of the once-punctured
torus. Hence, it induces a fibre-preserving involution
of any once-punctured torus bundle over the circle. The canonical decomposition
is preserved by this involution. We will prove the following
stronger version of Theorem 1.

\noindent {\bf Theorem 2.} {\sl Any ideal polyhedral
decomposition of a hyperbolic once-punctured torus bundle that 
is straight in the hyperbolic structure and that is invariant
under the fibre-preserving involution is equivariantly
isotopic to the monodromy ideal triangulation.}

\vskip 18pt
\centerline{\psfig{figure=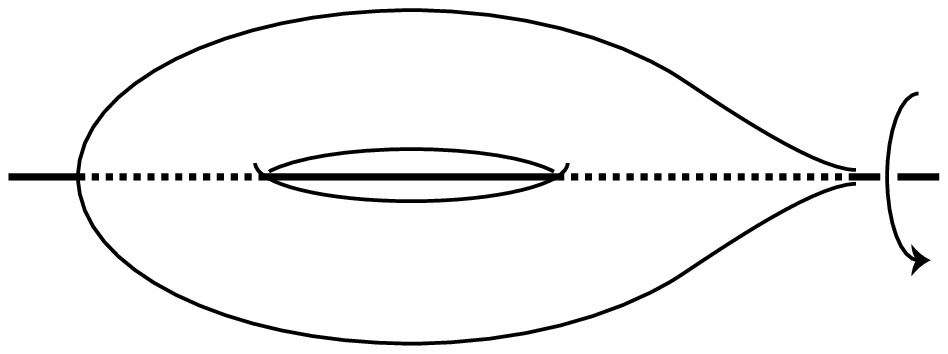,width=2in}}
\vskip 6pt
\centerline{Figure 1.}

The assumption that the decomposition is invariant under this
involution is vital. For example, the figure-eight knot 
complement, which is a once-punctured torus bundle, admits several straight ideal 
triangulations, as follows. Start with the canonical decomposition
into two ideal tetrahedra. Remove a common face, forming a single
ideal polyhedron. This can then be decomposed into three ideal tetrahedra
arranged around an edge, which form another straight ideal triangulation of the
knot complement.

A key concept in this paper is that of an
angled polyhedral decomposition of a 3-manifold $M$.
A {\sl polyhedron} is a 3-ball with a specified graph
in its boundary, which must satisfy the following conditions:
\item{$\bullet$} each vertex has valence at least three;
\item{$\bullet$} no edge is a loop;
\item{$\bullet$} each complementary region is a disc with
at least three edges in its boundary.

\noindent An {\sl ideal polyhedron} is a polyhedron with its
vertices removed. The resulting punctures are termed
{\sl ideal vertices}. An {\sl ideal polyhedral decomposition}
of a 3-manifold $M$ is a decomposition of $M - \partial M$ into ideal
polyhedra, with faces glued homeomorphically in pairs. An {\sl angled
polyhedral decomposition} is an ideal polyhedral decomposition,
together with an assignment to each edge of each ideal polyhedron
of an interior and exterior angle in the range $(0, \pi)$, which sum to
$\pi$. These must satisfy the following conditions:
\item{$\bullet$} the sum of the interior angles around
each edge is $2 \pi$;
\item{$\bullet$} for any closed curve in the boundary of
an ideal polyhedron that misses the ideal vertices,
that intersects each edge transversely at most once and that does
not lie wholly in a face, the sum of
the exterior angles of the edges it runs over
is at least $2 \pi$, with equality if and only it
it encircles an ideal vertex.

\noindent The second of the above conditions is, by Rivin's
theorem, exactly the requirement that an ideal polyhedron
with a given convex angle assignment to its edges
be realizable in hyperbolic space [7]. However, the hyperbolic
structures on these ideal polyhedra need not match up correctly
under the face identifications, and therefore an angled
polyhedral decomposition is significantly more general
than a decomposition of a hyperbolic 3-manifold into 
straight ideal polyhedra. A similar concept
(an angled spine) was introduced in [6].
We will in fact prove the following stronger version of
Theorem 2.

\noindent {\bf Theorem 3.} {\sl Any angled polyhedral
decomposition of a once-punctured torus bundle that is
invariant under the fibre-preserving involution is 
equivariantly ambient isotopic to the monodromy ideal 
triangulation.}

The strategy of the proof is to place the 1-skeleton of
$M$ in `thin position' in the fibration. We establish that 
either all the edges can be simultaneously equivariantly isotoped into fibres, or one 
of the fibres can be placed in a position that resembles almost normal 
form. We term this {\sl fairly normal form}, and it seems to
arise naturally when a 3-manifold admits certain sweep-outs
by surfaces with non-empty boundary. As in [6], it is
possible to assign an `area' to this 
fairly normal surface $F$, which is defined additively over the discs 
of intersection with the ideal polyhedra, and which turns out to be $-2 \pi \chi(F)$. 
However, using the fact that $F$ is invariant under the involution
and the fact that $F$ lies in a thick region of thin
position, we deduce that some disc area has more than 
$2 \pi$ or two discs have area more than $\pi$. 
It is simple to prove that each disc has non-negative
area. This gives a contradiction, since the Euler characteristic of 
a once-punctured torus is $-1$. We therefore
deduce that all the edges of the ideal polyhedral decomposition
can be made simultaneously level in the fibration.
With some further work, this implies that it is the monodromy
ideal triangulation. 

I would like to thank Brian Bowditch for introducing me to this
problem. I would also like to thank Makoto Sakuma who gave
a talk on his work with Akiyoshi, Wada and Yamashita [1] on 
quasi-fuchsian once-punctured torus groups, and the canonical
decomposition of two-bridge link complements. It seems likely
that the techniques of this paper can be applied to
two-bridge links.

\vskip 18pt
\centerline {\caps 2. The monodromy ideal triangulation}
\vskip 6pt

This section contains a description of the monodromy ideal triangulation
defined by Floyd and Hatcher [3].
The once-punctured torus $F$ can be constructed from two ideal
triangles by gluing their sides in pairs, and any ideal triangulation of
$F$ takes this form. Floyd and Hatcher found a very elegant way of
encoding the set of isotopy classes of such ideal triangulations as the
vertices of a tree. This tree is dual to a tessellation of the hyperbolic
plane by ideal triangles. The ideal vertices of this tessellation are
${\Bbb Q} \cup \{ \infty \}$ in the circle at infinity. Associated to each
ideal vertex, there is a properly embedded arc in $F$ with that
slope. Two ideal vertices are joined by a geodesic if and only if the
corresponding arcs can be isotoped off each other. These geodesics
divide the hyperbolic plane into ideal triangles, forming the
required tessellation, which is known as the {\sl diagram of $PSL(2,
{\Bbb Z})$}. (See Figure 2.) The ideal vertices of an ideal triangle correspond to
three disjoint non-parallel properly embedded
arcs in $F$, and hence an ideal triangulation. Thus, there
is one vertex of the dual tree for each isotopy class of ideal triangulation of the
once-punctured torus. Two vertices of the tree are joined
by an edge if and only if their corresponding ideal triangulations
differ by an {\sl elementary move}, which involves removing
one of the edges, resulting in a square with side identifications,
and then inserting the other diagonal of the square.

The monodromy of a once-punctured torus bundle induces a
homeomorphism of this tree. Any homeomorphism of a simplicial tree either has
a fixed point or leaves invariant a unique copy of ${\Bbb R}$,
known as the {\sl axis}. In the former case, the monodromy is
periodic, and hence the bundle is not hyperbolic.
In the latter case, we pick a vertex on the axis.
The unique path in the tree from this vertex
to its image under the monodromy homeomorphism
runs along the axis. It specifies a sequence of
elementary moves. These induce an ideal triangulation
of the bundle, as follows. Start with the once-punctured
torus with its initial ideal triangulation. Realize
the first elementary move by attaching an ideal
triangulation to the once-punctured torus, as
shown in Figure 3. Continue in this fashion
along the path in the graph as far as the final
vertex. We then glue the top and bottom ideal triangulations
via the monodromy. The result is the {\sl monodromy ideal 
triangulation} of the bundle.

\vskip 18pt
\centerline{\psfig{figure=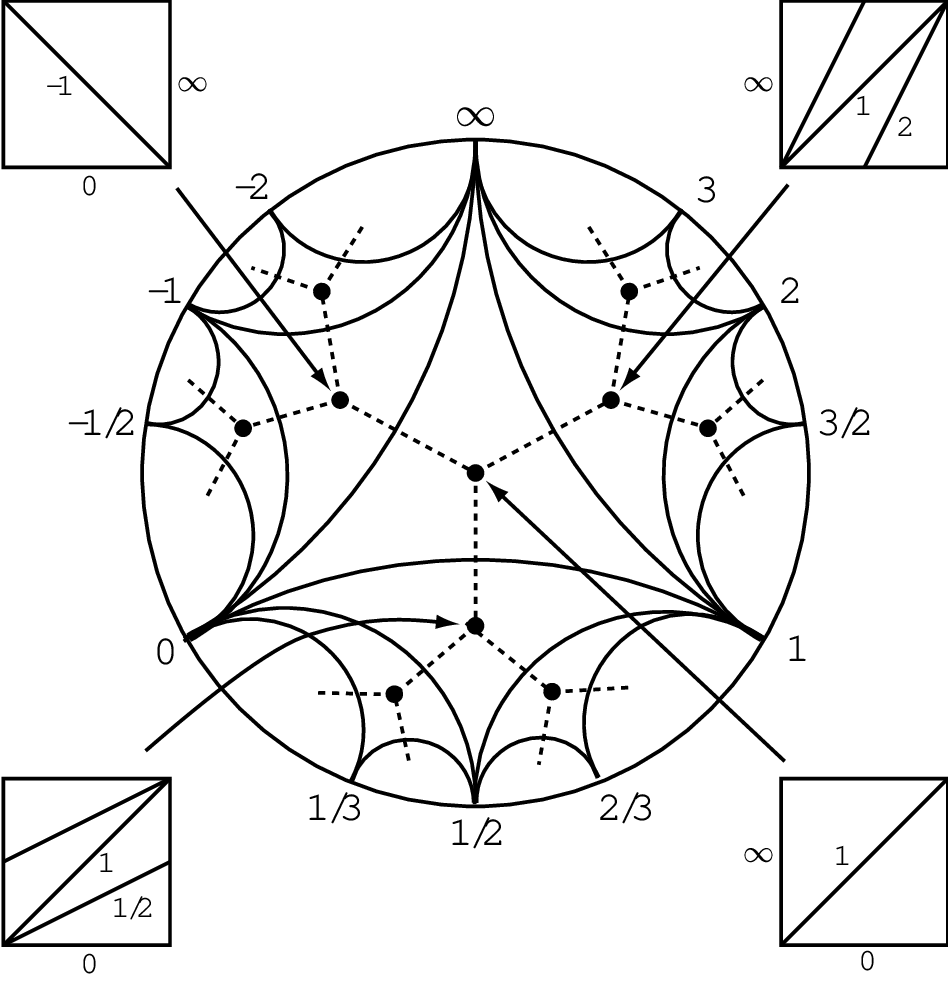,width=2.5in}}
\vskip 6pt
\centerline{Figure 2.}

\vskip 18pt
\centerline{\psfig{figure=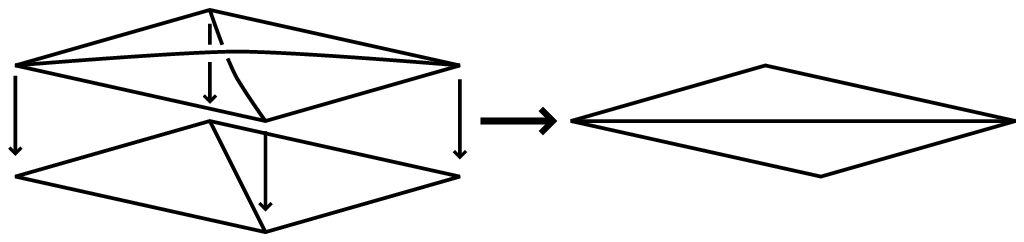}}
\vskip 6pt
\centerline{Figure 3.}

\vskip 18pt
\centerline{\caps 3. Fairly normal surfaces}
\vskip 6pt

It is a famous result, due to Rubinstein [8] and Stocking [9], that a 
strongly irreducible Heegaard surface of a compact orientable
3-manifold can be placed in almost normal form 
in any given triangulation. The arguments rely heavily on Gabai's concept of
thin position [4]. Many of these techniques generalise to other
settings: bridge decompositions of link complements, and surface fibrations
over the circle. However, it does not seem to be possible in general to deduce
that the relevant surface (be it a bridge punctured
2-sphere or a fibre) can be placed in almost normal form. But
there is a weaker version of the theory which we now
introduce.

Suppose we are given an ideal polyhedral decomposition
of a compact orientable 3-manifold with non-empty boundary. 
Truncate the ideal vertices, making the
polyhedra compact. These polyhedra now have two types of
faces: {\sl interior faces}, which are truncated copies
of the original faces, and {\sl boundary faces} which
are links of ideal vertices. The edges of the truncated
polyhedra also come in two types, {\sl boundary edges}
which lie in $\partial M$, and {\sl interior edges} which
we denote by $\Delta^1$.

A properly embedded surface $F$ is {\sl fairly normal}
if it intersects each polyhedron in fairly normal discs.
A {\sl fairly normal disc} is a properly embedded disc with the following properties:
\item{$\bullet$} it intersects any boundary face in
at most one arc;
\item{$\bullet$} it intersects any interior face in arcs,
that each starts and ends in distinct non-adjacent edges of the face;
\item{$\bullet$} it intersects any interior edge at most twice;
\item{$\bullet$} it is not parallel to a disc in the
boundary of the polyhedron disjoint from the interior
edges.

\noindent A {\sl weakly normal} disc is a properly embedded disc
satisfying the first two of the above conditions.
Note that in any given polyhedron, there are only
finitely many fairly normal disc types, but they are far more
numerous than normal discs. Some examples of fairly normal
discs in truncated tetrahedra are shown in Figure 4.

\vskip 18pt
\centerline{\psfig{figure=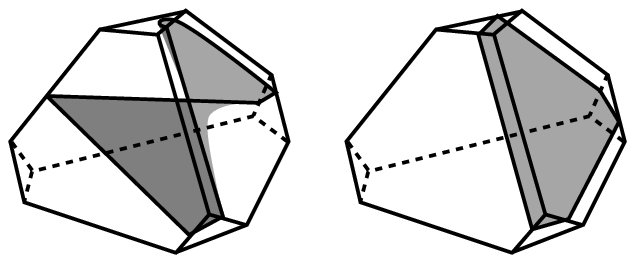}}
\vskip 6pt
\centerline{Figure 4.}

We will show that, when $M$ fibres over the circle,
and when the ideal polyhedral decomposition is angled,
then either there is an ambient isotopy taking each
interior edge into a fibre, or some fibre can be
placed in fairly normal form. This is nothing new: fairly 
normal surfaces are a generalisation 
of normal surfaces, and it is well known that a fibre can be ambient 
isotoped into normal form. However, the fairly normal surface
will have a number of
extra properties, which will eventually lead to a contradiction
in the case of a once-punctured torus.
For example, it will have the maximal number of intersection
points with $\Delta^1$ over all fibres not containing level edges.

When the ideal polyhedral decomposition is angled, it is
possible to assign an {\sl area} to weakly
normal (and more general) discs. This is the sum of the exterior angles of 
the interior edges the disc runs over (counted with multiplicity), 
plus the number of arcs
of intersection with the boundary faces multiplied
by $\pi$, then with $2 \pi$ subtracted. The {\sl
area} of a fairly normal surface $F$ is the sum of the areas of
its discs. It is shown in [6] that the area
of a properly embedded surface $F$ is equal to $-2 \pi \chi(F)$.

Many of the fairly normal surfaces $F$ we will consider
will have a {\sl face compression disc} which we
define to be a disc $D$ embedded in a polyhedron $P$, such that
\item{$\bullet$} $D$ is disjoint from $\Delta^1$;
\item{$\bullet$} the interior of $D$ is disjoint
from $\partial P \cup F$;
\item{$\bullet$} the boundary of $D$ is an arc in $F$,
an arc in an interior face and possibly an arc in a
boundary face;
\item{$\bullet$} the arc $\partial D \cap \partial P$
is not parallel in $\partial P - \Delta^1$ to a subarc
of $F \cap \partial P$.

\noindent For example, both the fairly normal discs shown in
Figure 4 have face compression discs.

If a weakly normal disc is surgered along a face compression
disc, the result is two discs that need not be weakly normal.
For they may have an arc of intersection with an interior face that has endpoints
\item{$\bullet$} in the same boundary edge,
\item{$\bullet$} in the same interior edge, or
\item{$\bullet$} in adjacent interior and boundary edges.

\noindent In the first case, the endpoints of this arc must close
up in the boundary face to form a closed curve that is disjoint from
the interior edges and that intersects the boundary faces in
only one arc. We call such a disc {\sl boundary-trivial}.

In the second and third cases, there is an obvious isotopy of
the discs that reduces the number of intersections with the
interior edges. We may repeat this procedure until we
end with discs that are boundary-trivial or weakly normal. Note
that in the third case, we always end with a weakly normal
disc.

We now define a certain type of fairly normal disc in a polyhedron $P$.
Let $\alpha$ be either an interior edge of $P$ or an arc properly
embedded in an interior face with endpoints in distinct boundary edges. 
Let ${\cal N}(\alpha)$ be a 
small regular neighbourhood of $\alpha$ in $\partial P$, and let 
$D$ be a properly embedded disc with boundary $\partial 
{\cal N}(\alpha)$. Then $D$ is a {\sl weak bigon}. It
is a {\sl bigon} if $\alpha$ was an interior edge. Note that
a fairly normal weak bigon is a bigon.

\noindent {\bf Lemma 4.} {\sl Let $D$ be a weakly normal disc.
Then the area of $D$ is non-negative. The area is zero if and only if it
is the link of an ideal vertex or a weak bigon.}

\noindent {\sl Proof.} 
Consider first the case where $D$ is disjoint from the boundary
faces. Suppose that it intersects some interior edge more
than once. It then has a face compression disc disjoint from
the boundary faces. If we
surger along this disc, the result is two new discs.
The area of $D$ is the sum of the areas of these two discs,
plus $2\pi$. If these discs fail to be weakly normal, there
is an isotopy that reduces their area and that takes
them to weakly normal discs. Thus, we may assume that $D$
intersects each interior edge at most once. But, then
by the definition of an angled polyhedral decomposition, its
area is non-negative and in fact is strictly positive
unless the disc is the link of an ideal vertex.

Suppose now that $D$ intersects the boundary faces
in a single arc. The interior edges emanating from the boundary
face containing this arc have exterior angles that sum to
$2 \pi$. Thus, there is a way of isotoping the arc off
the boundary face so that the new points of intersection
with the interior edges have total exterior angle at
most $\pi$. So, the resulting disc
has area at most that of $D$. It may fail to be
weakly normal, in which case there is an isotopy
which removes two intersection points with the interior
edges. Repeat until we end with a weakly normal disc $D'$.
By the above argument, $D'$ has non-negative area and hence so
does $D$. If the area of $D$ is zero, then $D'$ must have
been the link of an ideal vertex. Also, $D$ would have been
obtained from $D'$ by isotoping $D'$, without introducing
any intersection points with interior edges, so that a
single sub-arc of $\partial D$ runs over a boundary face.
The only way of doing this is to isotope $\partial D'$
along an interior edge or across an interior face. But in
both cases, the resulting $D$ has positive area.

If $D$ intersects the boundary
faces in more than two arcs, then the area is positive.
If it intersects the boundary faces exactly two times,
the area is non-negative, and is in fact strictly positive
unless it is disjoint from the interior edges. But in this
case, we claim that $D$ is a weak bigon. For, as above, there
is an isotopy of $D$ removing one of its arcs of intersections
with the boundary faces without increasing its area. After a further isotopy
that does not increase area, we end with a weakly normal or 
boundary-trivial disc. In the former case, this disc has positive
area since it has only a single arc of intersection with the
boundary faces. Hence $D$ also positive area. In the latter
case, $D$ was originally a weak bigon. $\square$

A consequence of the above lemma is that a fairly normal
surface cannot be a sphere or a properly embedded disc.
Also, by observing that bigons cannot be glued to links of
ideal vertices, we see that a fairly normal surface with
zero Euler characteristic is composed either entirely of
vertex links or entirely of bigons. Hence, it is a boundary-parallel
torus or a compressible annulus. This implies that a compact
orientable 3-manifold with an angled polyhedral decomposition
is irreducible, atoroidal, an-annular and has non-empty
boundary consisting of tori. Its interior therefore admits
a complete finite volume hyperbolic structure.

We now define some further types of fairly normal disc.
Let $\alpha$ be an arc in the boundary of a polyhedron $P$ such
that
\item{$\bullet$} the interior of $\alpha$ is disjoint from
the boundary faces;
\item{$\bullet$} the endpoints of $\alpha$ lie on 
boundary edges of distinct boundary faces;
\item{$\bullet$} the intersection of $\alpha$ with each
interior face is a collection of arcs, each of which
runs between distinct non-adjacent edges of the face;
\item{$\bullet$} $\alpha$ intersects each interior edge at most once.

\noindent Let ${\cal N}(\alpha)$ be a small regular neighbourhood
of $\alpha$ in $\partial P$, and let $D$ be a properly embedded
disc with boundary $\partial {\cal N}(\alpha)$. Then $D$
is an {\sl arclike} fairly normal disc, and $\alpha$ is 
its {\sl associated arc}. An arclike disc in a
truncated ideal tetrahedron is shown in Figure 5.

We give some further definitions. Let $D$ be a properly embedded 
disc parallel to a boundary face in a polyhedron $P$.
Pick one or two disjoint arcs $\alpha$ in $\partial P$,
starting in $\partial D$ and running across an interior face
to a boundary face. Ensure that each component of $\alpha$ is 
not parallel to a sub-arc of an interior 
edge, and that if $\alpha$ consists of two arcs, they are not
parallel. Modify by $D$ by isotoping along $\alpha$. 
The result is a {\sl modified vertex link}.

Similarly, let $D$ be two non-parallel ideal vertex links.
Let $\alpha$ be an arc in an interior face, running between the two discs of $D$.
Ensure that $\alpha$ is not parallel to a sub-arc of an interior 
edge. Modify $D$ by isotoping the boundaries of the two discs
towards each other along $\alpha$, and then fusing to form a single
disc, which we term a {\sl fused vertex link}.

\vskip 18pt
\centerline{\psfig{figure=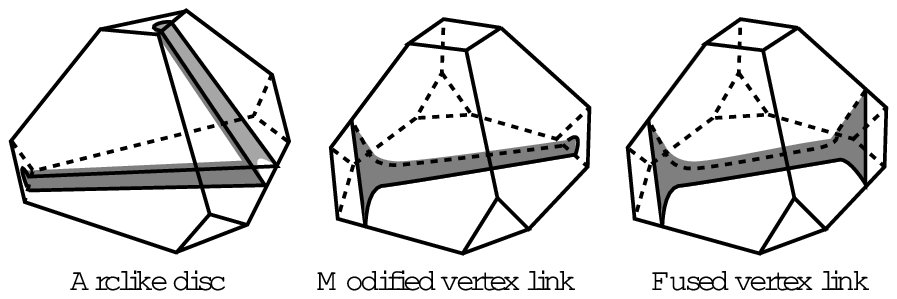}}
\vskip 6pt
\centerline{Figure 5.}

\noindent {\bf Lemma 5.} {\sl Let $D$ be a fairly normal disc in
a polyhedron $P$. Suppose that $D$ has a face
compression disc, and also that, if $D$ is arclike, the face
compression disc does not lie on the same side as its
associated arc. Then  the area of $D$ is at least $\pi$. Moreover, if $D$ 
is invariant under an involution of $P$ that preserves the disc's orientation and
transverse orientation, then it has area at least $2 \pi$.
If either of these inequalities are equalities, then the
disc must be disjoint from the interior edges or be a modified
vertex link or a fused vertex link.}

\noindent {\sl Proof.} Surger $D$ along the face compression disc to give two
discs $D_1$ and $D_2$ properly embedded in $P$. They need
not be weakly normal, but after an
isotopy that does not increase area, we end up with discs $D'_1$ and $D'_2$ that are 
weakly normal or boundary-trivial.

If $D'_1$ and $D'_2$ are both boundary-trivial, 
then $D$ was arclike and the face compression disc lay on the
arc side, contrary to hypothesis. So, at least one $D'_i$
is weakly normal, and so has non-negative area, by Lemma 4. If
either of $D'_1$ or $D'_2$ are boundary-trivial, then the face compression disc
was disjoint from $\partial M$. In this case, the
area of $D$ is equal to the area of $D_1$ plus the
area of $D_2$ plus $2 \pi$, and hence is at least $\pi$.
If both $D'_1$ and $D'_2$ are weakly normal, then they
have non-negative area. The area of $D$ is obtained by adding
the area of $D_1$ and $D_2$, and then adding either
$\pi$ or $2 \pi$, depending on whether the face compression
disc ran over a boundary face or not. Thus, we deduce that
the area of $D$ is at least $\pi$. 

Now suppose that this inequality is an equality. 
Then $D_1$ and $D_2$ must be equal to $D'_1$ and $D'_2$.
If $D_1$, say, is boundary-trivial, then $D_2$ must be a weakly 
normal disc with zero area. By Lemma 4, it is either
a vertex link or a weak bigon. Hence,
$D$ is either modified vertex link or disjoint from
$\Delta^1$. If $D_1$ and $D_2$ are both weakly normal,
then they must both have zero area, and the face compression
disc must have run over a boundary face. So, $D_1$ and $D_2$
must be weak bigons and hence, in this case, $D$ is disjoint from $\Delta^1$.

We must now show that if $D$ is invariant under an
involution that preserves its orientation and transverse orientation, 
then it has area at least $2 \pi$. The involution of $P$ must
be a rotation of order two, which preserves each component
of $P - D$. The involution takes the face compression disc to 
another. It is a straightforward matter to ensure that these
two discs are disjoint after an equivariant ambient isotopy.
If we now surger $D$ along both these discs, the result is three
discs $D_1$, $D_2$ and $D_3$, where $D_2$ is attached
to $D_1$ and $D_3$. So, $D_1$ and $D_3$ are swapped
by the involution and $D_2$ is invariant. There is an equivariant
isotopy that takes $D_1$ and $D_3$ to discs $D'_1$ and $D'_3$
that are weakly normal or boundary-trivial, and that takes
$D_2$ to a disc $D'_2$ that is weakly normal or that has
boundary in a single interior face. We call the
latter type of disc {\sl face-trivial}. Note that not
all the discs can be trivial, for then $D$ would be
arclike and its face compression disc would have been on the
arc side.

Consider first the case where the face compression discs are
disjoint from $\partial M$. Then the area of $D$ is the
sum of the areas of $D_1$, $D_2$ and $D_3$, plus $4 \pi$.
If $D'_2$ is face-trivial, its area is $-2\pi$ and the areas of $D'_1$ and $D'_3$
are non-negative. If $D'_1$ and $D'_3$ are boundary-trivial, then their
areas are $- \pi$, and the area of $D'_2$ is non-negative.
Thus, we deduce that the area of $D$ is at least $2 \pi$.
The case where the face compression discs intersect
$\partial M$ is easier. Each $D'_i$ is weakly normal, and
the area of $D$ is the sum of the areas of $D_1$, $D_2$
and $D_3$, plus $2 \pi$.
Hence, the area of $D$ is at least $2 \pi$.
The case of equality is straightforward and is omitted.
$\square$

\vskip 18pt
\centerline{\caps 4. Levelling the interior edges}
\vskip 6pt

In this section, we will prove the following result.

\noindent {\bf Theorem 6.} {\sl Let $M$ be a once-punctured torus
bundle with an angled polyhedral decomposition that is
invariant under the fibre-preserving involution.
Then there is an equivariant ambient isotopy taking each
interior edge into a fibre.}

\noindent {\sl Proof.} We place a transverse orientation on the
fibres, so that, locally, we may speak of one fibre being `above'
or `below' another.

We may equivariantly ambient
isotope the interior edges of $M$ so that
the following conditions hold. Each interior edge either lies entirely
in a fibre, or is transverse to the fibration at all but
finitely many critical points, which are local maxima or minima. 
We may assume that the endpoints of each non-level interior edge are critical.
We may also assume that a fibre cannot contain both level edges and
critical points, and that if a fibre contains more than one critical point
or more than one level edge, then
it contains precisely two and these are swapped by the involution.
The fibres containing critical points or level edges 
we term {\sl critical}. Critical fibres are divided into three types:
\item{$\bullet$} {\sl interior-critical},
which contain an isolated critical point in the interior of $M$;
\item{$\bullet$} {\sl boundary-critical},
which contain an endpoint of a non-level edge;
\item{$\bullet$} {\sl level-critical},
which contain a level edge.

\noindent Define the {\sl weight} of a critical fibre to be
\item{$\bullet$} two, if the fibre is interior-critical
and it contains two critical points (that must be swapped by the
involution);
\item{$\bullet$} the number of level edges it contains, if it
level-critical;
\item{$\bullet$} one, otherwise.

\noindent Let the {\sl width} of $\Delta^1$ be the sum, over all
critical fibres, of the product of the weight of the fibre
and its number of transverse points of intersection with
$\Delta^1$. Perform
an equivariant ambient isotopy which minimises width. 
The 1-skeleton is then in {\sl thin position}.

The reason for the weighting in the definition is as follows.
Suppose that local maxima of $\Delta^1$, in distinct adjacent
critical fibres, can be equivariantly 
isotoped past each other. Then we do not want this to alter
the width of $\Delta^1$. This is the case with this
choice of weighting. See Figure 6, for example.

\vskip 18pt
\centerline{\psfig{figure=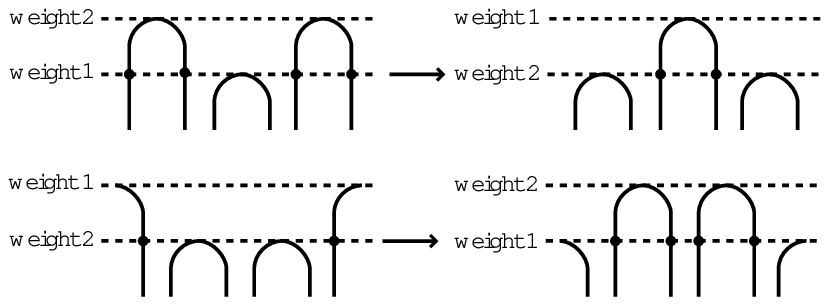}}
\vskip 6pt
\centerline{Figure 6.}

Suppose that not all the edges of $\Delta^1$ are level.
Consider an interval $I$ of fibres, disjoint from
the interior-critical and boundary-critical
fibres, starting just below a local maximum for $\Delta^1$,
and ending just above a local minimum for $\Delta^1$.
We may assume that each fibre in $I$, other than
the level-critical fibres,
has maximal intersection with $\Delta^1$ among all
fibres not containing level edges.
The interval $I$ is possibly divided up into sub-intervals
by level-critical fibres. Consider one such sub-interval $I'$,
which we initially take to be the highest sub-interval.
All fibres in this sub-interval are equivariantly ambient isotopic,
leaving $\Delta^1$ invariant. Hence, by examining the supremal
fibre, we see that each fibre in $I'$ has a strict upper disc, which is
defined as follows. An {\sl upper disc} (respectively,
{\sl lower disc}) for a fibre $F$ is an embedded disc $D$ in $M$ such that
\item{$\bullet$} its boundary is the union of an arc
in $F$, an arc in $\Delta^1$ and possibly an arc in $\partial M$;
\item{$\bullet$} the above arcs of $D \cap \Delta^1$ and $D \cap
\partial M$ are the only points of intersection between $D$ and
$\partial M \cup \Delta^1$;
\item{$\bullet$} near $\partial D \cap F$, its interior is
disjoint from $F$, and it emanates from the upper
(respectively, lower) side of $F$;
\item{$\bullet$} $D$ and its image under the fibre-preserving
involution are either equal or disjoint.

\noindent An upper or lower disc is {\sl strict} if its interior
is disjoint from $F$. Strict upper discs are shown in Figure 7. 

\vskip 18pt
\centerline{\psfig{figure=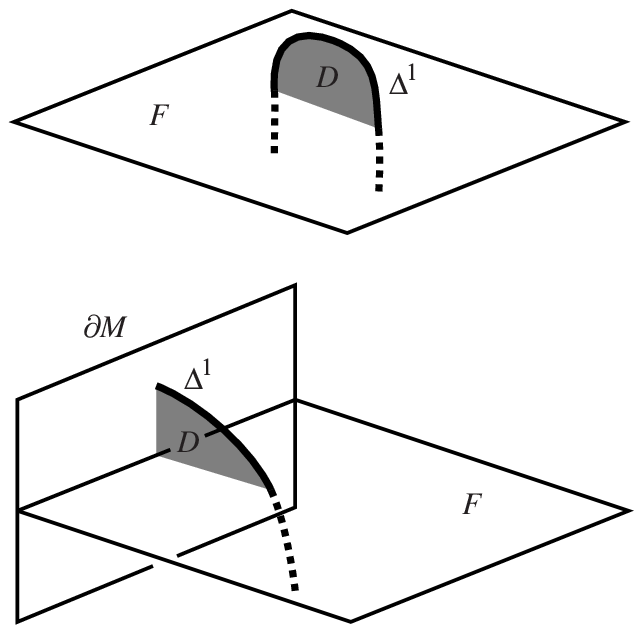}}
\vskip 18pt
\centerline{Figure 7.}

If, directly below $I'$, there is a level edge, then each fibre $F$ 
in $I'$ has a {\sl lower-parallelity disc}, which is an embedded disc
in $M$ such that
\item{$\bullet$} its boundary is an arc in $F$, an edge
of $\Delta^1$ and two arcs in $\partial M$;
\item{$\bullet$} its interior is disjoint from
$\Delta^1 \cup \partial M \cup F$;
\item{$\bullet$} it is attached to the lower side of $F$;
\item{$\bullet$} it and its image under the fibre-preserving
involution are either equal or disjoint.

\noindent {\sl Upper-parallelity} discs are defined similarly.

We define an upper disc and a
lower disc for a fibre $F$ to be a {\sl thinning
pair} if they are disjoint away from $F \cap \Delta^1$
and the same is true of the upper disc and
the image of the lower disc under the involution.
The following results are easy generalisations of standard
facts about thin position:
\item{$\bullet$} No fibre can have a thinning pair of
upper and lower discs.
\item{$\bullet$} No fibre can have an upper disc contained
in a lower disc or vice versa.
\item{$\bullet$} If a fibre contains a level edge of $\Delta^1$,
then it has no upper disc and no lower disc.

\noindent We will prove only the first of the above statements, as the
remainder have similar proofs. Suppose that $U$ and $L$ are upper and lower 
discs that form a thinning pair. Consider the intersection
between $U$ and the critical fibres, which
we may assume is a sub-arc of $\partial U$, a finite collection 
of points on $\partial U$, properly
embedded simple closed curves
and properly embedded arcs that are disjoint except possibly at their
endpoints. Let $D$ be a subdisc of $U$ separated off by an extrememost
arc $\alpha$. If the endpoints of $\alpha$ are equal, there
is an equivariant isotopy of $U$, keeping $\partial U$ fixed,
that removes $\alpha$. Otherwise, we equivariantly isotope
$D \cap \Delta^1$ across $D$ to remove $\alpha$, which does not increase the width of
$\Delta^1$. We can isotope $L$ similarly. Thus, we may
assume that $U$ and $L$ lie just above and just below a
critical fibre. There is then an equivariant ambient isotopy
that decreases the width of $\Delta^1$, either by moving
the maximum of $U$ past the minimum of $L$ or by cancelling
these critical points. This is impossible, since the 1-skeleton
is in thin position.

\noindent {\sl Claim 1.} In the interval $I'$, some fibre $F$ intersects each
interior face of $M$ transversely in the following components: simple closed curves,
arcs with endpoints in $\partial M$, arcs with endpoints
in distinct interior edges, and arcs with an endpoint in an interior
edge and a non-adjacent boundary edge.

We must show that some fibre $F$ in $I'$ intersects each face transversely
and has neither of the following components of intersection with
an interior face: an arc with endpoints in the same interior edge, or
an arc running from an interior edge to an adjacent boundary edge.
Note that, in each case, this arc separates off
an upper or lower disc. We term these {\sl local} upper and lower 
discs. It is impossible for any fibre to have
both local upper and local lower discs, by thin position.

During the isotopy specified by the interval $I'$, we may
assume that $F$ intersects each interior face transversely,
except at finitely points in $I'$, where the intersection
with any interior face performs the following moves:
\item{$\bullet$} add or remove an innermost simple
closed curve;
\item{$\bullet$} add or remove an extrememost arc with both its endpoints in
the same boundary edge;
\item{$\bullet$} move two curves together (possibly the same curve) until they
meet at a point, then resolve this singularity;
\item{$\bullet$} move a curve towards a boundary edge, until it
meets the edge at a point, then resolve this singularity;
\item{$\bullet$} move the endpoints of two arcs in the
same boundary edge towards each other until they meet,
and then resolve this singularity by pulling them away
from the boundary edge.

\noindent If the face is invariant under the involution, then two
copies of a move may occur simultaneously in the face.

The first two moves do not affect the existence of local
upper or lower discs. If two 
fibres differ by the third, fourth or fifth move, it is impossible 
for one to have a local lower disc and the other to have 
a local upper disc. For the non-transverse fibre between
them would have both such discs, and a small isotopy would
make these a thinning pair.

Now, the highest fibre of $I'$ has an upper disc with interior
disjoint from the interior faces. Hence, it cannot have a
local lower disc. Similarly, the lowest fibre 
of $I'$ cannot have a local upper disc. Hence, we deduce
that some fibre $F$ in $I'$ has no local lower and no local
upper disc, which proves the claim.

\noindent {\sl Claim 2.} There is an equivariant ambient
isotopy, that is fixed on the interior edges, taking $F$ into 
fairly normal form.

We will perform a series of moves that will each reduce
the number of components of intersection with the interior
faces. Hence, they will eventually terminate. 

Suppose that we were to perform equivariant compressions and boundary-compressions
to the components of intersection between $F$ and the polyhedra,
so that afterwards each such component is a disc that
intersects each boundary face in at most one arc.
The resulting surface we call $\overline F$. It is a collection
of spheres and discs, together with a single once-punctured torus.
The number of compressions and boundary-compressions used ($n$, say)
is equal to the number of spheres and discs. Each component 
of intersection between $\overline F$ and each polyhedron is
either weakly normal, boundary-trivial or face-trivial.
The area of each sphere and disc is negative, and
so, by Lemma 4, it must contain a boundary-trivial or
face-trivial disc. But a face-trivial or boundary-trivial 
disc can only be attached to another such disc.
So, there are in total $2n$ trivial discs in
$\overline F$. The original surface $F$ is obtained from 
$\overline F$ by attaching $n$ tubes, at least one of
which must be attached to the once-punctured torus component.
Hence, if $n>0$, some trivial disc lies in $F$. We can perform an equivariant isotopy 
of $F$ to remove this disc. The resulting surface still satisfies
Claim 1. Thus, we can repeat until $F$ is composed of weakly normal
discs. If one of these discs fails to satisfy the last condition in
the definition of being fairly normal, there is
an obvious equivariant isotopy in the complement of
the interior edges that reduces the number of intersections
with the boundary faces. Finally, if a
disc of $F$ intersects an interior edge more than twice, then
a thinning pair of upper and lower discs can be found,
contradicting thin position. This proves the claim.

\noindent {\sl Claim 3.} If a fairly normal fibre $F$ has a strict upper disc,
then it has a face compression disc on the upper side of $F$.
If $F$ has an upper-parallelity disc, then
either there is a face compression disc on the upper side
of $F$, or this upper-parallelity disc
can be ambient isotoped into a polyhedron.
Similar statements are true for discs below $F$.

Consider the intersection of the strict upper disc $U$ with the interior faces. 
We may assume that, near $\Delta^1$, the only intersection between
$U$ and the interior faces is $\Delta^1 \cap U$. Thus,
the intersection is $\Delta^1 \cap U$ and a collection of arcs
and simple closed curves properly embedded in $U$, disjoint from $\Delta^1$.
The closed curves may be removed. Consider an extrememost arc in
$U$ away from $\Delta^1 \cap U$. If it has both endpoints in
$\partial M$, it may be removed. The arc separates off a disc that 
lies in a polyhedron $P$. If this is not
a face compression disc, then its intersection with $\partial P$
is parallel in the complement of $\Delta^1$ to an
arc in $F \cap \partial P$. There is then an ambient
isotopy of $U$, removing this arc of intersection. Repeating 
this process, we isotope $U$ into a polyhedron $P$.
A further small isotopy in $P$ makes it a face compression disc. 

The situation with an upper-parallelity disc is 
similar. The intersection between this disc and the interior
faces is again a collection of simple closed curves and arcs, and by the above
argument, we may assume that it consists only of arcs running from $\partial M$
to $\partial M$. We may successively remove the arc
closest to the edge of $\Delta^1$, until the parallelity
disc lies entirely in a single polyhedron. This proves the
claim.

So far we have not used the fact that $F$ is a once-punctured
torus. We will do so now, by means of an area argument. 
We know that the fairly normal fibre $F$ in $I'$
has a strict upper disc. Hence,
by Claim 3, it has a face compression disc on its upper side.
Suppose first that this disc does not emanate from the
arc side of an arclike disc. Then, by Lemma 5,
some fairly normal disc of $F$ has area at least $\pi$. 
Moreover, if it invariant under the involution, then
it has area at least $2 \pi$. This disc and its image
under the involution therefore account for all the area
of $F$. By Lemma 5, each is disjoint from $\Delta^1$,
or is a modified or fused vertex link. The remaining
discs have zero area, and hence by Lemma 4, each is a
bigon or vertex link.

\noindent {\sl Case 1.} The discs with positive area
are disjoint from $\Delta^1$.

A vertex link disc cannot be attached to a disc disjoint
from $\Delta^1$. Therefore, there are no vertex link discs.
So, $F$ is entirely disjoint from $\Delta^1$. But
$F$ was, by construction, a fibre with maximal intersection 
with $\Delta^1$. This is a contradiction, which proves the
theorem in this case.

\noindent {\sl Case 2.} The discs with positive area
are modified vertex links.

A modified vertex link has arcs of intersection with interior
faces that run from an interior edge to a boundary face.
These arcs cannot be attached to a vertex link or to
a bigon. So, this arc must be attached to a modified
vertex link on the other side of the face. But this 
creates a component of $\partial F$ which bounds a
disc in $\partial M$, which does not occur. Thus, this
case does not arise.

\noindent {\sl Case 3.} The only disc with positive area
is a fused vertex link.

A bigon cannot be attached to a fused vertex link or to 
a vertex link. So, $F$ consists only of a fused vertex link
and vertex links. This implies that $F$ is closed, which is
a contradiction.

Thus, we may assume that $F$ has an arclike disc with the
arc on the upper side. It cannot therefore have an arclike
disc with the arc on the lower side. For this would imply
that it has a thinning pair of lower and upper discs, contradicting
thin position.

Now, either $F$ has a strict lower disc or there
are level edges directly below $F$ in the fibration. In the former case,
the above argument proves the theorem. So suppose that directly
below $F$ there are level edges. These give rise to
lower-parallelity discs $P$. If, below $F$, there is a face
compression disc, we are done. Hence, by Claim 3,
we may assume that the lower-parallelity discs can each be isotoped
into a polyhedron. 

\noindent {\sl Claim 4.} $F$ has at most two arclike discs, and if
it has exactly two, then these are swapped by the involution.
Also, the lower-parallelity discs $P$ are attached to these arclike discs. 

If not, then we can find an arclike disc $D$ disjoint from $P$. This arclike
disc gives rise to an upper disc for the fibre containing the level
edges in $P$. This contradicts thin position, and so establishes the claim.

The fibre $F'$ directly below the level edges is obtained
by isotoping $F$ across $P$. By Claim 4, $F'$ is fairly normal
and has no arclike discs. If, directly below $F'$, there is a
local minimum for $\Delta^1$, then we may apply the
above argument to $F'$ and hence prove the theorem.
If below $F'$ there are more level edges, then this
gives rise to further lower-parallelity discs. We may assume
that these can be isotoped into a polyhedron, for
otherwise the theorem is proved. But this implies that the
parallelity discs above $F'$ are disjoint from those below
$F'$. They therefore extend to parallelity discs below $F$
that are disjoint from its arclike discs. Without
changing the width, we may equivariantly isotope the level
edges below $F'$ and the level edges above $F'$ past each other.
This creates a fibre containing a level edge and having an
upper disc, which contradicts thin position.
$\square$

\vfill\eject
\centerline{\caps 5. Conclusion of the proof}
\vskip 6pt

\noindent {\bf Lemma 7.} {\sl Every interior edge is invariant under the
involution, and its orientation is reversed.}

\noindent {\sl Proof.} Each interior edge $e$ lies in a fibre which is preserved
under the involution. This involution leaves invariant every properly
embedded essential arc in the fibre, up to isotopy, and reverses its orientation.
Thus, $e$ and its image either coincide or are
parallel in the surface. However, it is impossible for distinct
edges in a fibre to be parallel. For, the parallelity disc between
two adjacent edges lies in the complement of $\Delta^1$. The
intersection between this disc and the interior faces is a collection
of properly embedded arcs and simple closed curves. By removing
innermost curves and extrememost arcs, the disc may be isotoped
into a polyhedron. We then see that $e$ must coincide with its
image. $\square$

\noindent {\bf Proposition 8.} {\sl The fixed point set of the involution
cannot have an arc of intersection with any face.}

\noindent {\sl Proof.} Via standard Morse theory, we may apply
an equivariant isotopy of the interior faces, keeping the interior edges
fixed, so that the interior of each face has only finitely many critical points
in the fibration. This gives a singular
foliation on each face. We define the {\sl critical leaves}
to be those that contain critical points. So, the complement
of the critical leaves has a product foliation. In the interior 
of the faces, we may take the critical points to have the following standard
forms: 
\item{$\bullet$} saddles, which form a 4-valent vertex of a 
critical leaf;
\item{$\bullet$} maxima and minima, whose corresponding critical leaf
is a point.

\vskip 18pt
\centerline{\psfig{figure=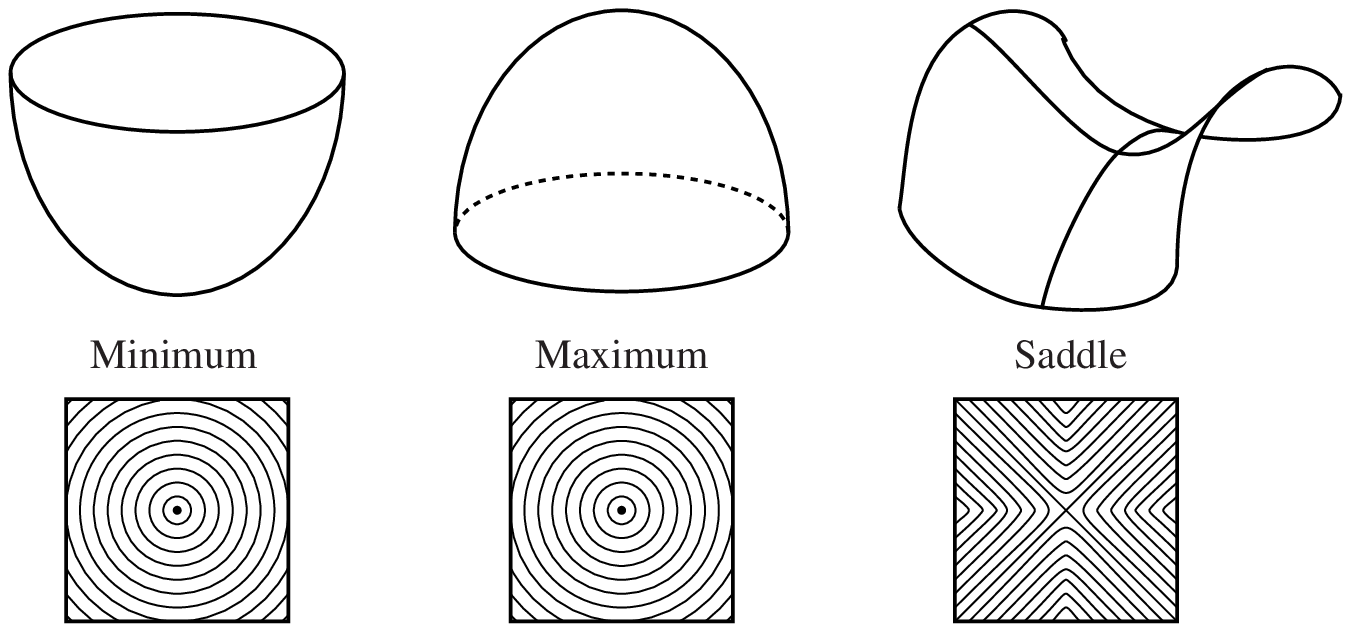,width=3in}}
\vskip 18pt
\centerline{Figure 8.}

\noindent We insist that
each boundary edge is either level or has finitely many critical
points of the following forms:
\item{$\bullet$} maxima and minima, which create a critical leaf
that is a point;
\item{$\bullet$} half-saddles, which contribute 2-valent vertices
to a critical leaf.

\vskip 18pt
\centerline{\psfig{figure=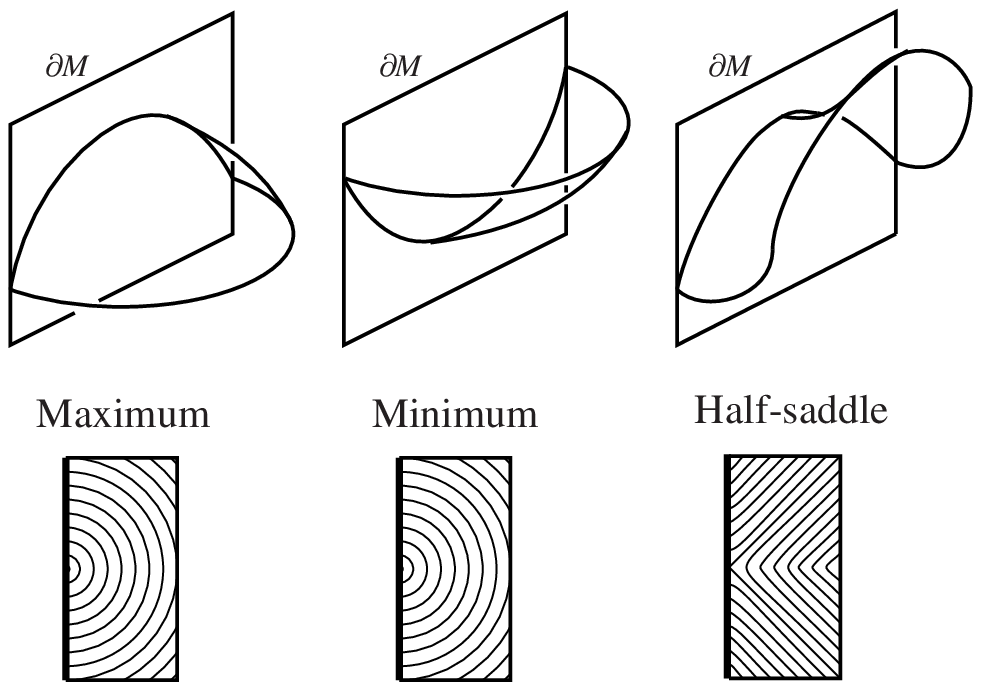,width=2.5in}}
\vskip 18pt
\centerline{Figure 9.}

Near each level edge,
we may arrange that each interior face has a finite number of {\sl switches}, 
shown in Figure 10. We may assume that each fibre not containing level edges
contains at most two critical points, and that if it does contain
two, then these are swapped by the involution. We may also
arrange that fibres containing level edges are disjoint from
the saddle, half-saddle, maximum and minimum critical points.

\vskip 18pt
\centerline{\psfig{figure=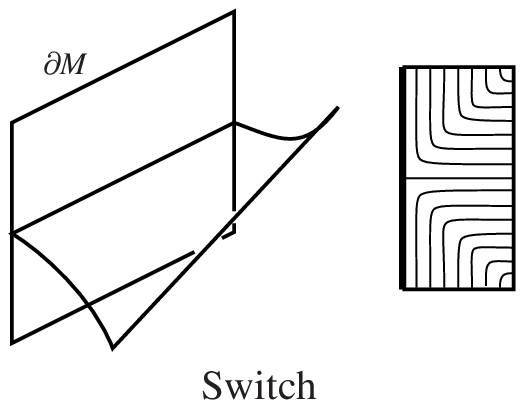,width=1.5in}}
\vskip 6pt
\centerline{Figure 10.}

We may remove all maxima and minima as follows, at the cost of
introducing {\sl ridges}. These are properly embedded arcs in an
interior face, with endpoints in $\partial M$, which are level
in the fibration and form either the highest or lowest points in
a neighbourhood of the arc.

\vskip 18pt
\centerline{\psfig{figure=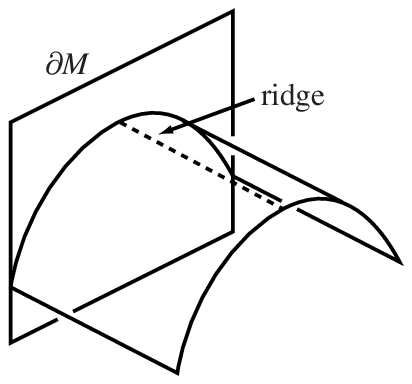,width=1.2in}}
\vskip 6pt
\centerline{Figure 11.}

For suppose that there is a local maximum, say. The adjacent
complementary region $R$ of the critical leaves has a product
foliation consisting of properly embedded arcs or circles,
depending on whether the maximum lies in the boundary or
the interior of $M$. Hence, $\partial R$ contains no other maxima
or minima. So, it contains another type of critical
point. There are a number of possibilities
for the critical point (or points) in $\partial R$ and for the
maximum. In each case, we can cancel the maximum with a critical point 
(or points). For example, if the maximum lies on the boundary
of $M$, and the other critical point is a half-saddle,
then these may be cancelled to form a ridge. The only
small complication is when the maximum is in the interior
of the face and there is a single critical point in $\partial R$
that is a saddle. It may be the case that $\partial R$ runs
over this saddle twice. But then $R$ separates
off a subdisc of the face, which contains at least one
maximum or minimum. So, by passing to an innermost region,
we can assume this case does not arise.
Thus, we may ensure that each face now has no maxima or minima. This implies
that the complement of the critical leaves is foliated by
arcs, and that each complementary region has two critical
leaves in its boundary.

Suppose now that a face contains an arc of intersection
with the fixed point set of the involution. Then, the face is subjected to a
rotation about this arc. Since each interior edge is reversed by the involution, 
the gluing pattern of the edges of the face
must be as shown in Figure 12, possibly with further
identifications. 
The arc of the fixed point set has endpoints in two interior edges of the face, the
`top' and `bottom' edges, $t$ and $b$. Consider the complementary region $R$ of
the critical leaves attached to the part of $t$ that intersects
the fixed point set, and consider
the critical leaf in $\partial R$ not intersecting $t$. Its
intersection with the fixed point set is not a critical point,
since the fixed point set is transverse to the fibration.
Hence, emanating from this point is a sub-arc $\alpha$ of the
critical leaf. There are the following possibilities
for the next critical points along this arc:
\item{1.} $\alpha$ is properly embedded in the face and ends in switches.
\item{2.} $\alpha$ ends in saddles.
\item{3.} $\alpha$ ends in half-saddles.
\item{4.} $\alpha$ contains a subset of the bottom edge $b$.

\vskip 18pt
\centerline{\psfig{figure=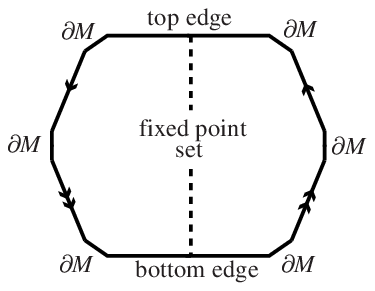}}
\vskip 6pt
\centerline{Figure 12.}

\vskip 18pt
\centerline{\psfig{figure=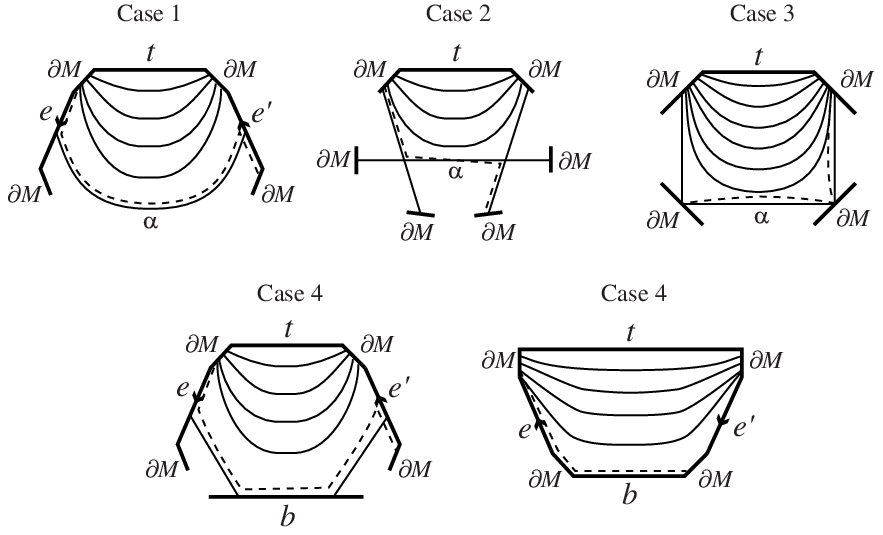,width=4in}}
\vskip 6pt
\centerline{Figure 13.}

Let us consider Case 1. The switches lie in interior edges $e$ and $e'$,
neither of which can be $t$ or $b$. These must be
the same edge in $M$, since each interior edge is preserved under the
involution. The arc $\alpha$ and the edge $e = e'$ both lie in the same fibre $F$.
Consider the arc $\beta$ in the face that starts in $\partial M$, runs 
forwards along $e$ as far as an endpoint of $\alpha$, then runs along $\alpha$,
then runs backwards along $e'$ back to $\partial M$. This is
shown in Case 1 of Figure 13 as a dotted arc. Since $e$
is an essential embedded arc in $F$ and the interior of $\alpha$ is
embedded, $\beta$ is homotopic in
$F$ to a properly embedded arc. Hence, it and its
image under the involution are homotopic in $F$ as unoriented
arcs. Thus, we deduce that
$\beta$ and its image are homotopic in $M$ via a homotopy
whose interior misses $\Delta^1$. However, this is impossible, by the following
claim, since $\beta$ and its image are distinct essential arcs in the
face, not parallel to an interior edge. The proof of this claim is simple and omitted.

\noindent {\sl Claim.} Suppose that two properly embedded
arcs in interior faces are homotopic, via a homotopy that
misses $\Delta^1$. Then either these are parallel arcs in
the same face, or they are both parallel in the faces
to interior edges, or they are both inessential in their
faces.

Cases 2 and 3 are dealt with similarly. In each case,
we apply the argument of Case 1 to the dotted arcs
in Figure 13. Consider now Case 4. The non-singular leaves in $R$ end in boundary edges.
These cannot be adjacent to both $t$ and $b$. Suppose that
they are not adjacent to $b$. Then, emanating from $b$ in
the boundary of $R$, there are either switches or level
boundary edges. In the former case, the critical leaf
continues from the switches to other interior edges $e$
and $e'$, which are swapped by the involution.
We can apply the argument of Case 1 to the dotted arc
shown in Figure 13. In the latter case, the level boundary edges
are attached to other interior edges $e$ and $e'$.
To apply the argument of Case 1 to the dotted arc in Figure 13, we
need to know that it is homotopic in its fibre $F$ to a properly embedded
arc. This is clear if $e$ and $b$ are distinct interior
edges. We claim that they cannot be the same edge.
Suppose that they were. If the orientations of $e'$, $b$
and $e$ were all consistent around $\partial R$, then
their concatenation would not be a primitive element
of $H_1(F,\partial F)$, and hence would not be homotopic
to a properly embedded arc. However, their concatenation
is homotopic to such an arc, as can be seen by considering
a nearby non-singular leaf. If the orientations of $e$
and $b$ are inconsistent around $\partial R$, then the
level edge between them would have to start and end at
the same point, and hence would close up to form a level
closed curve in $\partial M$. But, it and its
image under the involution have disjoint interiors,
which would then not be possible. So, we can again
apply the argument of Case 1. $\square$

Denote the number of polyhedra by $p$, the number of interior faces by $f$
and the number of interior edges by $e$. When the fixed point set is
cut along the points of intersection with the
2-skeleton, the result is $s$ arcs, say. We have the following
inequalities:
\item{$\bullet$} $p \geq s$, since each polyhedron can contain
at most one arc of the fixed point set in its interior. This is
an equality if and only if each polyhedron is invariant under
the involution.
\item{$\bullet$} $s \geq e$, as Lemma 7 implies that each interior
edge intersects the fixed point set, and Proposition 8 gives that
no face contains an arc of intersection with the fixed point set.
This is an equality if and only if the fixed point set is disjoint from the interior
of the faces.
\item{$\bullet$} $e \geq p$, as each polyhedron has
at least four faces, and each interior face lies on
the boundary of two polyhedra (or forms two faces of
a single polyhedron). Hence, $2f \geq 4p$. But by Euler
characteristic, $e - f + p = 0$ and so this gives the
required inequality. This is an equality if and only
if each ideal polyhedron is an ideal tetrahedron.

We then deduce that the above inequalities must be
equalities. Therefore, each polyhedron is a tetrahedron
that is preserved by the involution. The fixed point set
must run through every tetrahedron exactly once
between the midpoints of opposite edges. The interior faces
of each tetrahedron are partitioned into two orbits, and
the two faces in an orbit form a square with side identifications.
The involution applies an order two rotation to this
square, and reverses the orientation of each of the interior
edges. Hence, the side identifications of the square are
those of a once-punctured torus. We therefore see that the ideal triangulation
is obtained from a copy of the once-punctured torus by successively
attaching ideal tetrahedra, realizing elementary moves, and then gluing
top to bottom. In this sequence of moves, no move can be immediately
followed by its inverse, as this would create an edge with valence
two, which is impossible in angled polyhedral decomposition.
Hence, this is the monodromy ideal
triangulation of a representation of $M$ as a once-punctured torus
bundle. However, $M$ fibres over the circle in only one way,
since it easy to check that $H_2(M, \partial M) \cong {\Bbb Z}$.
As the space of ideal triangulations of the punctured torus is
a tree, there is only one possible such sequence of moves that
realizes the monodromy.  
Thus, this is the monodromy ideal triangulation of $M$ for the given
bundle structure. This proves Theorems 1, 2 and 3.

\vfill\eject
\centerline {\caps References}
\vskip 6pt

\item{1.} {\caps H. Akiyoshi, M. Sakuma, M. Wada and Y. Yamashita},
{\sl Ford domains of punctured torus groups and two-bridge knot
groups}, in Knot Theory, Conference Proceedings, Toronto 1999.
\item{2.} {\caps D. B. A. Epstein and R. C. Penner}, 
{\sl Euclidean decomposition
of non-compact hyperbolic manifolds}, J. Differential Geom. {\bf 27} (1988) 
67-80.
\item{3.} {\caps W. Floyd and A. Hatcher}, {\sl Incompressible
surfaces in punctured torus bundles}, Topology Appl. {\bf 13}
(1982) 263--282.
\item{4.} {\caps D. Gabai}, {\sl Foliations and the topology of
3-manifolds III}, J. Differential Geom. {\bf 26} (1987)
479--536.
\item{5.} {\caps T. J{\o}rgensen}, {\sl On pairs of once-punctured
tori}, unfinished manuscript.
\item{6.} {\caps M. Lackenby}, {\sl Word hyperbolic Dehn surgery},
Invent. Math. {\bf 140} (2000) 243--282.
\item{7.} {\caps I. Rivin}, {\sl On geometry of convex ideal polyhedra in
hyperbolic 3-space}, Topology {\bf 32} (1993) 87--92.
\item{8.} {\caps J. H. Rubinstein}, {\sl Polyhedral minimal surfaces,
Heegaard splittings and decision problems for 3-dimensional
manifolds}, Proceedings of the Georgia Topology Conference,
AMS/IP Stud. Adv. Math, vol. 21, Amer. Math Soc. (1997) 1--20.
\item{9.} {\caps M. Stocking}, {\sl Almost normal surfaces in 3-manifolds},
Trans. Amer. Math. Soc. {\bf 352} (2000) 171--207.
\item{10.} {\caps A. Thompson}, {\sl Thin position and the recognition
problem for the 3-sphere}, Math. Res. Lett. {\bf 1} (1994) 613--630.
\item{11.} {\caps J. Weeks}, {\sl SnapPea},  available at
http://thames.northnet.org/weeks/index/

\vskip 18pt

\+ Mathematical Institute, Oxford University, \cr
\+ 24-29 St Giles', Oxford OX1 3LB, England. \cr

\end